\documentclass[11pt,reqno]{amsart}

\usepackage[text={160mm,240mm},centering]{geometry}            
\usepackage{tikz}
\usepackage{tikz-3dplot}
\usetikzlibrary{arrows.meta}
\usetikzlibrary{commutative-diagrams}
\usetikzlibrary{cd}
\usepackage{amssymb}
\usepackage{amsmath}
\usepackage{longtable}
\usepackage{textcomp}
\usepackage{xcolor}
\definecolor{navyblue}{RGB}{25,45,80}
\usepackage[colorlinks=true,
            linkcolor=navyblue,
            citecolor=navyblue,
            urlcolor=navyblue]{hyperref}
\newtheorem{theorem}{Theorem}[section]
\newtheorem{lemma}[theorem]{Lemma}
\newtheorem{proposition}[theorem]{Proposition}

\newtheorem{corollary}[theorem]{Corollary}

\theoremstyle{definition}
\newtheorem{definition}[theorem]{Definition}

\newtheorem*{claim}{Claim}

\theoremstyle{remark}
\newtheorem{remark}[theorem]{Remark}

\numberwithin{equation}{section}

\usepackage{enumitem}
\setlist[enumerate,1]{label=(\roman*)}

\newcommand{\m}{\mathfrak{m}}
\newcommand{\n}{\mathfrak{n}}
\newcommand{\p}{\mathfrak{p}}
\newcommand{\chara}{\mathrm{char}}
\newcommand{\HH}{\mathcal{H}}
\newcommand{\wcl}{\mathrm{w.cl}}
\newcommand{\cl}{\mathrm{cl}}
\DeclareMathOperator{\coker}{coker}
\DeclareMathOperator{\im}{im}

\begin{document}

\title{Growth of $h$-Functions and $F$-Thresholds}

\author{Hongrui Ma}
\address{Department of Mathematical Sciences,
Tsinghua University, Beijing, 100084, P. R. China.}
\email{mhr23@mails.tsinghua.edu.cn}

\author{Huaiqing Zuo}
\address{Department of Mathematical Sciences,
Tsinghua University,
Beijing, 100084, P. R. China.}
\email{hqzuo@mail.tsinghua.edu.cn}





\begin{abstract}
We establish a quadratic lower bound for the increments of the $h$-function $h_s(I,J)$ on suitable intervals, where $I$ and $J$ are $\mathfrak m$-primary ideals in a Noetherian local ring of positive characteristic and dimension at least two. This complements the upper bound given by Lipschitz continuity and yields general lower bounds for $F$-thresholds. These bounds are strict in regular local rings of dimension at least three. In particular, if $R$ is
regular of dimension $n\geq 3$ and $\mu_R(I)\leq n+1$, then $c^I(I)>(n+1)/2$. Applications include results for $s$-closures, strict inequalities relating colength and multiplicity, the quotient of Milnor number by Tjurina number of isolated hypersurface singularities, and bounds involving $F$-pure and log canonical thresholds.

\par\smallskip
\begingroup
\makeatletter
\renewcommand{\keywordsname}{Keywords}
\renewcommand{\subjclassname}{MSC 2020}
\noindent\@setkeywords\par
\noindent\@setsubjclass\par
\makeatother
\endgroup

\end{abstract}

\maketitle

\section{Introduction}
Hilbert–Samuel multiplicity is a fundamental numerical invariant in commutative algebra that measures the complexity of singularities. In positive characteristic, the Frobenius endomorphism provides another approach to studying singularities, giving rise to Hilbert–Kunz multiplicity.

Throughout, rings are commutative and Noetherian. More precisely, let \((R,\m)\) be a Noetherian local ring of dimension \(n\), let \(I\) be an \(\m\)-primary ideal, and let \(M\) be a nonzero finitely generated \(R\)-module. The Hilbert-Samuel multiplicity $e_R(I;M)$ of \(M\) with respect to \(I\) (or $e(I;M)$ when $R$ is understood, $e(I)$ when $M=R$) is defined by $$e_R(I;M)=\lim_{t \rightarrow \infty} \frac{n!\lambda_R(M/I^tM)}{t^n}.$$
If \(R\) has positive characteristic \(p>0\), the Hilbert–Kunz multiplicity $e_{HK}(I;M)$ of \(M\) with respect to \(I\) (or $e_{HK}(I)$ when $M=R$) is defined by
\[
e_{HK}(I;M):=
\lim_{e\to\infty}
\frac{\lambda_R\!\left(M/I^{[p^e]}M\right)}{p^{en}},
\]
where \(I^{[p^e]}=(a^{p^e}\mid a\in I)\).

The theory of $h$-functions, or $s$-multiplicities, first introduced in \cite{e-and-eHK}, connects Hilbert–-Samuel and Hilbert–Kunz multiplicities and provides, to some extent, a unified framework for their study. Meng and Mukhopadhyay \cite{Meng} developed the theory and obtained interesting results. The $h$-function is defined as follows:
\begin{definition}\label{def-h}\textup{(\cite[Definition 2.5]{e-and-eHK})}
    Let $(R,\mathfrak{m})$ be a local ring of dimension $n$ and characteristic $p>0$. $I,J \subset R$ are $\mathfrak{m}$-primary ideals. Let $M$ be a finitely generated $R$-module. For $s>0$, the limit
    $$h_s(I,J;M)=\lim_{e \rightarrow +\infty} \frac{\lambda_R(M/(I^{\lceil sp^e \rceil}+J^{[p^e]})M)}{p^{en}} $$ exists.
    We write $h_s(I,J)$ when $M=R$.
\end{definition}
Taylor proved $h_s(I,J)$ is Lipschitz continuous \cite[Theorem 2.7]{e-and-eHK}. That is, there exists a constant $C>0$ such that for every $b>a>0,\ h_b(I,J)-h_a(I,J) \leq C|b-a|.$ In our Section \ref{section-quadratic-l-bound}, we prove that the increment of $h_s$ has a quadratic lower bound on a suitable interval. More precisely, let $(R,\mathfrak{m},k)$ be a local ring of dimension $n \geq 2$ and characteristic $p>0$ and let $I,J$ be two $\m$-primary ideals, $J^M \subseteq I$ for some integer $M \geq 1$. Set $$\HH_s(n)=\frac{1}{n!}\sum_{i=0}^{\lfloor s \rfloor}(-1)^i\binom{n}{i}(s-i)^n$$ for $s\geq 0$ and $\HH_s(n)=0$ for $s \leq 0$. Set $$v=\frac{1}{M},\ r=\mu_R(I)+\mu_R((I+J)/I)-n.$$ If $I$ has a reduction $Q$ that is a parameter ideal, we can furthermore choose $r=\mu_R((J+Q)/Q)$. Set $$\Phi_{s,v,r}(n)=\HH_{s}(n)-r\HH_{s-v}(n).$$ With this notation, our first main result is as follows.

\begin{theorem}\label{intro-thm-h-l-bound}\textup{(Theorem \ref{thm-ha-hb-l-bound}, Corollary \ref{cor-hb-ha-qua})}
    For $0<a<b\leq v+\frac{n}{2}$, we have
    $$h_b(I,J)-h_a(I,J) \geq e_R(I)(\Phi_{b,v,r}(n)-\Phi_{a,v,r}(n)).$$

    In particular, if $r \geq 1$, let $\tau$ be the unique solution in the interval $(v,n)$ of $$\Phi_{s,v,r}'(n)=\HH_{s}'(n)-r\HH_{s-v}'(n)=0,$$ where $\Phi_{s,v,r}'(n)$ denotes the derivative of \(\Phi_{s,v,r}(n)\) with respect to \(s\). Then $0<v<\tau\leq\frac{n+v}{2}<v+\frac{n}{2}$. For any fixed real number \(0<\delta<\tau\), there exists a constant $C>0$ such that $$h_{b}(I,J)-h_a(I,J) \geq C(b-a)^2$$  for every $\delta \leq a \leq b \leq \tau$.
\end{theorem}

The quadratic lower bound restricts where the \(h\)-function can become constant, and therefore yields a lower bound for the corresponding \(F\)-threshold.

\begin{theorem}\textup{(Theorem \ref{thm-cJI}, Proposition \ref{prop-condition-hold}, Corollary \ref{cor-cII})}
    Retain the notation of Theorem \ref{intro-thm-h-l-bound}, assume $n \geq 3,\ r \geq 1$, and suppose that one of the following conditions holds:
    \begin{enumerate}
        \item\label{item-condition-regular} $R$ is regular.
        \item $R$ is Cohen-Macaulay and $J$ is a parameter ideal.
    \end{enumerate}
    Then $$c^J(I)>\tau.$$
    In particular, we have
    \begin{equation}\label{cJI-ineq-eq}
        c^J(I) > \begin{cases}
            \dfrac{n+1}{2M}, &\mu_R(J) \leq n+1,\\
            \dfrac{\mu_R(J)^{1/(n-1)}}{M(\mu_R(J)^{1/(n-1)}-1)},&\mu_R(J)\ge n+2.
        \end{cases}
    \end{equation}
    As a corollary, if $R$ is regular and $Q$ is a $\m$-primary ideal with $\mu_R(Q)\leq n+1$, then $$c^Q(Q)>\frac{n+1}{2}.$$
\end{theorem}

The second application concerns the $s$-closure of an ideal $I$, denoted by $I^{\cl_s}$ for $s \geq 1$, see \cite{e-and-eHK}. For an isolated hypersurface singularity $f$ over a field of characteristic zero with its Jacobian ideal $j(f)$, Zhang and Zuo \cite[Corollary 6.2]{HNS-conj} proved $$j(f):f \nsubseteq (f,j(f)).$$ They also proved that $$j(f):f \nsubseteq \overline{j(f)}$$ may not hold, where $\overline{j(f)}$ denotes the integral closure of $j(f)$ \cite[Theorem 7.3]{HNS-conj}. In this paper, we prove

\begin{theorem}\label{intro-thm-nsubset-s-closure}\textup{(Theorem \ref{thm-nsubset-s-closure})}
    Let $(R,\m)$ be a complete regular local ring of dimension $n \geq 3$ and characteristic $p>0$. Let $I$ be an $\m$-primary ideal of $R$, with a reduction $Q$ that is a parameter ideal. Let $r=\mu_R(I)-n \geq 1.$ Let $1<\tau \leq \frac{n+1}{2}$ be the unique solution of $\Phi'_{s,1,r}(n)=\HH_s'(n)-r\HH_{s-1}'(n)=0$. Then $Q:I \nsubseteq Q^{\cl_{n-\tau}}$.
\end{theorem}
As a corollary, we get
\begin{corollary}\textup{(Corollary \ref{cor-mu-s-closure})}
    Let $f \in k[[x_1,\dots,x_n]]$ define an isolated hypersurface singularity, $n \geq 3$, $\chara k>0$. Let $j(f)$ denote its Jacobian ideal. Assume that $f\in \overline{j(f)}$. Then $j(f):f \nsubseteq j(f)^{\mathrm{cl}_{\frac{n-1}{2}}}$.
\end{corollary}
\begin{remark}
    This confirms our conjecture \cite[Conjecture 3.15]{mu-tau}.
\end{remark}

Since $j(f) \subseteq (j(f))^* \subseteq j(f)^{\cl_s} \subseteq \overline{j(f)}$ for $s \geq 1$, where $j(f)^*$ denotes the tight closure of $j(f)$, we provide a more precise result on how large $(j(f):f)$ is. This also completely answers the questions concerning the quotient of the Milnor number to the Tjurina
number for hypersurface singularities in arbitrary characteristic, which is partly done in \cite{mu-tau}.
\begin{theorem}\textup{(Corollary \ref{cor-mu-tau-p})}
    For an isolated hypersurface singularity $f$ of dimension $n-1$, $n \geq 3$, we have
    $$\frac{\mu(\mathcal{O}_f)}{\tau(f)}<\frac{1}{\mathcal{H}_{\frac{n+1}{2}}(n)-\mathcal{H}_{\frac{n-1}{2}}(n)},$$
    and if $f$ is defined over a field of characteristic zero, then 
    $$\frac{\mu(f)}{\tau(f)}<\frac{1}{\mathcal{H}_{\frac{n+1}{2}}(n)-\mathcal{H}_{\frac{n-1}{2}}(n)},$$
    where $\mu(f),\mu(\mathcal{O}_f),\tau(f)$ denote the Milnor number, generalized Milnor number and Tjurina number.
\end{theorem}
\begin{remark}
    The characteristic-zero case and the corresponding non-strict inequality in positive characteristic were established in \cite[Theorem 1.4]{mu-tau}. Here we strengthen the latter to a strict inequality for $n\ge 3$.
\end{remark}

The same argument also yields a strict inequality between colength and multiplicity in regular local rings, strengthening the bounds in \cite[Proposition 5.7]{Lech} (or \cite[Theorem 2.10]{lower-bound-s}, \cite[Theorem 3.1]{e_hk-eI-equality}).

\begin{corollary}\textup{(Corollary \ref{cor-colength-e-p}, Corollary \ref{cor-colength-e-0})}
    Let $(R,\m,k)$ be an equicharacteristic regular Noetherian local ring of dimension $n \geq 3$. Let $I$ be an $\m$-primary ideal of $R$ such that $r=\mu_R(I)-n \geq 1.$ Then we have $$\lambda_R(R/I)> e(I)(\HH_s(n)-r\HH_{s-1}(n)) $$ for $s \geq 1$.
\end{corollary}

As a corollary, we get
\begin{corollary}\textup{(Corollary \ref{cor-F-regular}, Corollary \ref{cor-klt})}
    Let \((R,\mathfrak m)\) be an \(F\)-finite regular local ring of characteristic \(p>0\) (resp. a regular local ring essentially of finite type over a field of
characteristic zero) and dimension \(n=3\). Let $I$ be an $\m$-primary ideal. If $\mu_R(I) \leq 4$ and the pair $(R,I^1)$ is not strongly $F$-regular (resp. not Kawamata log terminal), then $\lambda_R(R/I) \geq 19$.
\end{corollary}
This shows that strengthening the inequality to a strict one is both necessary and nontrivial.

\section{Ideal quotients, duality and length}\label{section-lemma}
Let $(R,\m,k)$ be a Noetherian local ring with maximal ideal $\m$ and residue field $k$. Let $E_R(k)$ denote the injective hull of $k$. For an $R$-module $M$, denote by $M^\vee=\mathrm{Hom}_R(M,E_R(k))$ the Matlis dual of $M$. Denote $\lambda_R(M)$ the length of $M$, $\mu_R(M)$ the minimal number of generators of $M$. Note that the Matlis duality preserves length, that is, $\lambda_R(M)=\lambda_R(M^\vee)$ for every module $M$ of finite length.

We begin with several properties of the quotient of ideals and their duality.

\begin{lemma}\label{lemma-quotient}
    Let $(R,\m,k)$ be a Noetherian local ring with maximal ideal $\m$ and residue field $k$. Let $I$ be an ideal of $R$.
    \begin{enumerate}
    \item\label{item-quotient} Let $J$ be an ideal of $R$. Then $(I:J)/I=0:_{R/I} J$. If $K \subseteq I$ is an ideal, then $(I:J)/K=(I/K):_{R/K}((J+K)/K)$.
    \item\label{M/IM-duality} Let $M$ be an $R$-module of finite length. Then $(M/IM)^\vee \cong 0:_{M^\vee} I,\ (0:_M I)^\vee \cong M^\vee/IM^\vee.$ In particular, for $c \in R$, we have $\lambda_R(0:_M c)=\lambda_R(M/cM)$.
    \end{enumerate}
\end{lemma}

\begin{proof}
    (i) It can be easily checked by definition.\\
    (ii) $(M/IM)^\vee=\{\varphi \in M^\vee:\ \varphi(IM)=0\}$, while
    \begin{equation}
        \begin{aligned}
            \varphi(IM)=0 &\Leftrightarrow \varphi(az)=0,\ \forall a \in I, z \in M \\
            & \Leftrightarrow a\varphi=0, \forall a \in I \\
            & \Leftrightarrow \varphi \in 0:_{M^\vee}I.
        \end{aligned}
    \end{equation}
    
    This gives $(M/IM)^\vee \cong 0:_{M^\vee} I$. The second isomorphism follows from $M^{\vee\vee} \cong M$, \cite[Proposition 3.2.12(c)]{cohen-macaulay-ring}. Putting $I=(c)$, we have $\lambda_R(0:_M c)=\lambda_R(M^\vee)-\lambda_R(cM^\vee)=\lambda_R(M)-\lambda_R(cM)=\lambda_R(M/cM)$.
\end{proof}

\begin{lemma}\label{lemma-length-eq}
    Let $(A,\n)$ be a Noetherian local domain of dimension $d \geq 1$.  Let $\varphi: M \rightarrow N$ be a morphism of finite $A$-modules. Moreover, there exists an integer $m \geq 1$, such that $\n^mM=0,\ \n^mN=0$. Assume that there exists $c \in A$ such that $c \neq 0,\ c$ is not a unit and $c \ker \varphi=c\coker \varphi=0$. Then $$|\lambda_A(M)-\lambda_A(N)| \leq \max\{\mu_A(M),\mu_A(N)\}\lambda_A(\frac{A}{\n^m+(c)}).$$ In particular, if $A,c$ are fixed, $M,N$ vary with $m$ and $\mu_A(M),\ \mu_A(N)$ are bounded independently of $m$, then $\lambda_A(M)-\lambda_A(N)=O(m^{d-1})$.
\end{lemma}

\begin{proof}
    Since $\n^mM=\n^mN=0$, $M,N$ are of finite length.

    Since $c \ker \varphi=0$, $\ker \varphi \subseteq 0:_M c.$ Thus $$\lambda_A(\ker \varphi) \leq \lambda_A(0:_M c)=\lambda_A(M/cM)$$ by Lemma \ref{lemma-quotient}\ref{M/IM-duality}. Since $M$ is generated by $\mu_A(M)$ elements and $\n^mM=0$, we have a surjection $$(A/(\n^m+(c)))^{\mu_A(M)} \twoheadrightarrow M/cM.$$ Therefore, $$\lambda_A(\ker \varphi) \leq \lambda_A(M/cM) \leq \mu_A(M) \lambda_A(A/(\n^m+(c))).$$

    On the other hand, $\coker \varphi$ is a quotient of $N$, and is killed by $c$. Thus $N \twoheadrightarrow \coker \varphi$ factors through $N/cN$. Similarly we have the surjection $$(A/(\n^m+(c)))^{\mu_A(N)} \twoheadrightarrow N/cN \twoheadrightarrow \coker \varphi,$$ which gives $$\lambda_A(\coker \varphi) \leq \mu_A(N) \lambda_A(A/(\n^m+(c))).$$ By the exact sequence 
    $$0 \rightarrow \ker \varphi \rightarrow M \rightarrow N \rightarrow \coker \varphi \rightarrow 0,$$ we get $$|\lambda_A(M)-\lambda_A(N)|=|\lambda_A(\ker \varphi)-\lambda_A(\coker \varphi)| \leq \max\{\mu_A(M),\mu_A(N)\}\lambda_A(\frac{A}{\n^m+(c)}).$$

    Since $A$ is a domain, $c \neq 0, \ c$ is not a unit, $\dim A/(c) \leq d-1$. We have $$\lambda_A(A/(\n^m+(c)))=\lambda_A(\frac{A/(c)}{\n^m(A/(c))})=O(m^{d-1}).$$ The last result follows.
\end{proof}

Let $(R,\m,k)$ be a Noetherian local ring. Let $J \subset R$ be an $\m$-primary ideal of $R$. Then $B=R/J$ is an Artin local ring. Therefore, the canonical module $\omega_B=(R/J)^\vee=E_B(k)$. Denote $\m_B=\m/J$ the maximal ideal of $B$.

\begin{lemma}\label{lemma-duality}
    Let $(R,\m,k)$ be a Noetherian local ring, let $J,B,\m_B$ be as above.
    \begin{enumerate}
        \item\label{item-canonical} $\omega_B$ is annihilated by $J$. $\mu_R(\omega_B)=\dim_k \frac{\omega_B}{\m \omega_B}=\dim_k \frac{J:\m}{J}=\dim_k 0:_{R/J} \m$, where $0:_M \m$ is denoted by $\mathrm{Soc}(M)$, $\dim_k 0:_M \m$ is denoted by $r(M)$ for an Artinian $R$-module $M$.
        \item\label{item-gorenstein} $B$ is Artinian Gorenstein if and only if $\dim_k 0:_B\m_B=r(B)=1$.
        \item\label{item-0:0:I} If $R$ is regular, $J$ is a parameter ideal, then for any ideal $I \subseteq R$, $J:(J:I)=J+I$.
    \end{enumerate}
\end{lemma}

\begin{proof}
    (i) Since $\omega_B$ is an $R/J$-module, $J\omega_B=0$. The formula for $\mu_R(\omega_B)$ follows from \cite[Proposition 3.2.12(d)]{cohen-macaulay-ring}. For the notations $r(M)$ and $\mathrm{Soc}(M)$, we also follow \cite{cohen-macaulay-ring}.\\
    (ii) It follows from \cite[Theorem 3.2.10]{cohen-macaulay-ring}.\\
    (iii) In this case, $B$ is a complete intersection, hence Artinian Gorenstein by \cite[Proposition 3.1.20]{cohen-macaulay-ring}. Then by \cite[Lemma 1.11(1)]{0:0:I}, $0:_B(0:_B \overline{I})=\overline{I}$ in $B$, where $\overline{I}=(I+J)/J$ is the image of $I$ in $B$. It follows that $J:(J:I)=J+I$ in $R$.
\end{proof}

The next lemma uses Matlis duality to convert an inclusion \(I^{a+b}\subseteq J\) into a length estimate involving the canonical module of \(R/J\).

\begin{lemma}\label{lemma-long-ineq}
    Let $(R,\m,k)$ be a Noetherian local ring. Let $I \subset R$ be an ideal, $J \subset R$ be an $\m$-primary ideal. Denote $\omega_{R/J}$ the canonical module of $R/J$. If there exist integers $a,b \geq 0$ such that $I^{a+b} \subseteq J$, then
    \begin{equation}\label{lemma-ineq-1}
        \lambda_R(\frac{I^a+J}{J})\leq \lambda_R(\frac{\omega_{R/J}}{I^b\omega_{R/J}})\leq \mu_R(\omega_{R/J}) \lambda_R(\frac{R}{I^b+J}).
    \end{equation}
    In particular, 
    \begin{equation}\label{lemma-ineq-2}
        \begin{aligned}
            \lambda_R(R/J) &\leq \lambda_R(\frac{R}{I^a+J})+\lambda_R(\frac{\omega_{R/J}}{I^b\omega_{R/J}})\\
            & \leq \lambda_R(\frac{R}{I^a+J})+\mu_R(\omega_{R/J}) \lambda_R(\frac{R}{I^b+J}).
        \end{aligned}
    \end{equation}

\end{lemma}

\begin{proof}
    By $I^{a+b} \subseteq J$, we have $I^a \subseteq J:I^b$. Therefore, $$(I^a+J)/J \subseteq (J:I^b)/J=0:_{R/J} I^b$$ by Lemma \ref{lemma-quotient}\ref{item-quotient}. Taking lengths, we get 
    \begin{equation}\label{first-ineq-1}
        \lambda_R((I^a+J)/J) \leq \lambda_R(0:_{R/J} I^b).
    \end{equation}
    By Lemma \ref{lemma-quotient}\ref{M/IM-duality}, $$(0:_{R/J} I^b)^\vee \cong \frac{(R/J)^\vee}{I^b(R/J)^\vee}=\frac{\omega_{R/J}}{I^b\omega_{R/J}}.$$ Since the Matlis duality keeps the length, we have $$\lambda_R(0:_{R/J} I^b)=\lambda_R(\frac{\omega_{R/J}}{I^b\omega_{R/J}}).$$ Combining this with (\ref{first-ineq-1}), we get the first inequality in (\ref{lemma-ineq-1}).

    Since $\omega_{R/J}$ has $\mu_R(\omega_{R/J})$ generators and $J\omega_{R/J}=0$ by Lemma \ref{lemma-duality}\ref{item-canonical}, we have a surjection $(R/J)^{\mu_R(\omega_{R/J})} \twoheadrightarrow \omega_{R/J}$, which induces a surjection $$(\frac{R}{I^b+J})^{\mu_R(\omega_{R/J})} \twoheadrightarrow \frac{\omega_{R/J}}{I^b\omega_{R/J}}.$$ This gives the second inequality in (\ref{lemma-ineq-1}). (\ref{lemma-ineq-2}) then follows easily.
\end{proof}

\section{A quadratic lower bound for increments of $h$-functions}\label{section-quadratic-l-bound}

Let $(R,\m)$ be a Noetherian local ring of dimension $n$ with characteristic $p>0$. Let $I,J$ be two $\m$-primary ideals. The theory of $h$-functions (Definition \ref{def-h}) provides a connection between Hilbert–Samuel multiplicity, Hilbert–Kunz multiplicity, and $F$-thresholds.

Denote $\nu^J_I(p^e)=\max\{m:I^m \nsubseteq J^{[p^e]}\}$, $\mu^J_I(p^e)=\min\{m:J^{[p^e]} \nsubseteq I^m\}$. There exist limits $$c^J(I)=\lim_{e \rightarrow \infty} \frac{\nu^J_I(p^e)}{p^e} \qquad \mathrm{and} \qquad b^J(I)=\lim_{e \rightarrow \infty} \frac{\mu^J_I(p^e)}{p^e}$$ that follow from \cite[Theorem 3.4]{cJI}, \cite[Lemma 3.2]{e-and-eHK}.
We have the following properties.

\begin{proposition}\textup{(\cite[Proposition 2.6, Theorem 2.7, Lemma 3.3]{e-and-eHK})}\label{prop-h}
    Let $(R,\mathfrak{m})$ be a local ring of dimension $n$ and characteristic $p>0$. $I,J \subset R$ are $\mathfrak{m}$-primary ideals. Let $M$ be a finitely generated $R$-module.
    \begin{enumerate}
        \item The function $h_s(I,J;M)$ is Lipschitz continuous.
        \item If $s \leq b^J(I)$ then $h_s(I,J;M)=\frac{s^n}{n!}e(I;M)$. If $s \geq c^J(I)$ then $h_s(I,J;M)=e_{HK}(J;M)$.
        \item\label{item-h-integral} If $I_1$ and $I_2$ have the same integral closure, then $h_s(I_1,J;M)=h_s(I_2,J;M)$.
    \end{enumerate}
\end{proposition}

The Lipschitz continuity provides a linear upper bound for the increments of $h_s$. In this section, we establish a quadratic lower bound for these increments on a suitable interval. 

We first record a faithfully flat extension that allows us to assume completeness and an infinite residue field while preserving the numerical invariants used below.

\begin{lemma}\label{ring-extension}
    Let $(R,\m,k)$ be a Noetherian local ring of dimension $n$. Let $R_1=R[t]_{\m R[t]}$, $S=\hat{R_1}$, the $\m R_1$-adic completion of $R_1$. Then 
    \begin{enumerate}
        \item $(S,\m S)$ is a complete Noetherian local ring with infinite residue field.
        \item $R \rightarrow S$ is faithfully flat. 
        \item For a finite length $R$-module $W$ and a finite $R$-module $L$, $\lambda_S(W \otimes_R S)=\lambda_R(W),\ \mu_S(L \otimes_R S)=\mu_R(L)$.
        \item For an $\m$-primary ideal $I \subset R$, $e_S(IS)=e_R(I)$. In particular, $\dim S=n$.
        \item If \(R\) has characteristic \(p>0\), then for $\m$-primary ideals $I,J \subseteq R,$ $h_s(I,J;R)=h_s(IS,JS;S)$.
    \end{enumerate}
\end{lemma}

\begin{proof}
    (i) Obvious.\\
    (ii),(iii),(iv) $$\frac{L \otimes_R S}{\m S(L \otimes_R S)}=(L/\m L)\otimes_k k(t)$$ gives $\mu_S(L \otimes_R S)=\mu_R(L)$. The remaining assertions follow from \cite[Lemma 8.4.2(2),(4),(6)]{Integral-closure-ideal} and the fact that completion preserves flatness, lengths and multiplicities.\\
    (v) Frobenius powers commute with extension of ideals, hence we have $(JS)^{[q]}=J^{[q]}S$ for $q=p^e$. Therefore, $$\lambda_R(\frac{R}{I^{\lceil sq \rceil}+J^{[q]}})=\lambda_S(\frac{S}{(IS)^{\lceil sq \rceil}+(JS)^{[q]}}).$$ The result follows.
\end{proof}

Then we establish a structure lemma.
\begin{lemma}\label{lemma-rankAR}
    Let $(R,\m,k)$ be a complete Noetherian local ring of dimension $n$ with characteristic $p>0$. Let $Q=(x_1,\dots,x_n)$ be an $\m$-primary ideal of $R$. Fix a coefficient field \(k\subseteq R\), and set $A=k[[x_1,\dots,x_n]] \subset R$, $\n=(x_1,\dots,x_n)$ the maximal ideal of $A$. Then $R$ is a finite $A$-module, $\n R=Q$ and $\lambda_R(M)=\lambda_A(M)$ for every $R$-module $M$ of finite length. Moreover, $\mathrm{rank}_A R=e_R(Q)$, denoted by $\rho$. Then there exists an injection $A^\rho \hookrightarrow R$ and $0 \neq c \in A$ such that $cR \subseteq A^{\rho}$.
\end{lemma}
\begin{proof}
    $R$ is a finite $A$-module by Cohen Structure Theorem \cite[Theorem 29.4]{Matsumura-commutative}. $\lambda_R(M)=\lambda_A(M)$ since $R$ and $A$ have the same residue field. $\mathrm{rank}_A R=e_R(Q)$ follows from \cite[Corollary 11.2.6]{Integral-closure-ideal}.

    Let $u_1,\dots,u_\rho \in R$ form a basis of $R \otimes_A \mathrm{Frac}A$. The map
    $$A^\rho \hookrightarrow R,\ (a_1,\dots,a_\rho) \mapsto a_1u_1+\dots+a_\rho u_\rho $$ is injective (since $u_1,\dots,u_\rho$ are linearly independent over $\mathrm{Frac}A$). In the following, we regard $A^\rho$ as a submodule of $R$. Then $(R/A^\rho) \otimes_A \mathrm{Frac}A=0$. Therefore, every generator of $R/A^\rho$ is killed by some nonzero elements in $A$. Thus, there exists $0 \neq c \in A$ such that $c(R/A^\rho)=0.$ That is, $cR \subseteq A^\rho$.
\end{proof}

The following property allows us to enlarge the second ideal when estimating increments from below.

\begin{lemma}\label{lemma-hb-ha-J}
    Let $(R,\mathfrak{m})$ be a local ring of dimension $n$ and characteristic $p>0$. $I,J_1,J_2 \subset R$ are $\mathfrak{m}$-primary ideals, $J_1 \subseteq J_2$. Then for $0<a<b$, we have
    $$h_b(I,J_1)-h_a(I,J_1) \geq h_b(I,J_2)-h_a(I,J_2).$$
\end{lemma}
\begin{proof}
    For $q=p^e$, $J_1^{[q]} \subseteq J_2^{[q]}$ induces a surjection
    $$\frac{I^{\lceil aq\rceil}+J_1^{[q]}}{I^{\lceil bq\rceil}+J_1^{[q]}} \twoheadrightarrow \frac{I^{\lceil aq\rceil}+J_2^{[q]}}{I^{\lceil bq\rceil}+J_2^{[q]}}.$$
    To see this, note that every element of the right-hand side can be written as $z+\delta,\ z \in I^{\lceil aq \rceil},\ \delta \in J_2^{[q]}$, which has a preimage on the left hand side. Taking lengths, dividing both sides by \(q^n\) and passing to the limit, we obtain the result.
\end{proof}

\begin{lemma}\label{lemma-J-Q-kq}
    Let $(R,\mathfrak{m})$ be a local ring of dimension $n$ and characteristic $p>0$ and let $I,J$ be two $\m$-primary ideals, $J^M \subseteq I$ for some integer $M \geq 1$. Suppose $Q$ is a minimal reduction of $I$. Then there exists $w \in \mathbb{N}$ such that $I^{w+1}=QI^w$. 
    
    Set $r=\mu_R((J+Q)/Q)$. Then there exist $y_1,\dots,y_r \in J$ such that $J+Q=Q+(y_1,\dots,y_r)$. Set $q=p^e,\ k_q=\lfloor q/M \rfloor-w.$ For sufficiently large $q$, we have $k_q \geq 0,\ y_i^q \in Q^{k_q}$.
\end{lemma}
\begin{proof}
    $I^{w+1}=QI^w$ follows by \cite[Definition 1.2.1]{Integral-closure-ideal}. Choosing a group of minimal generators of $(J+Q)/Q$ and lifting to $J$, we get $y_1,\dots,y_r \in J$ such that $J+Q=Q+(y_1,\dots,y_r)$. Write $q=Mu+t,\ u=\lfloor\frac{q}{M}\rfloor.$ When $q$ is sufficiently large, $u$ is large enough so that $k_q=u-w \geq 0$. We have $y_i^q=(y_i^M)^u \cdot y_i^t \in I^u \subseteq Q^{u-w}=Q^{k_q}$.
\end{proof}

We next establish the following dimension inequality.
\begin{lemma}\label{lemma-dim-ineq}
    Let $k$ be an arbitrary field. $n \geq 1,\ q \geq 2$ are integers. Denote $$M_q=\frac{k[[x_1,\dots,x_n]]}{(x_1^q,\dots,x_n^q)}.$$ Regard $M_q$ as a $k$-vector space generated by the monomials $x^{\alpha}$, where $\mathbf{x}^{\alpha}=x_1^{\alpha_1}\dots x_n^{\alpha_n}$ for $\alpha=(\alpha_1,\dots,\alpha_n) \in \mathbb{N}^n$. Denote $$V_{q,j}=\mathrm{span}\{\mathbf{x}^{\alpha}:0\leq \alpha_i<q,\ \alpha_1+\dots+\alpha_n \leq j\} \subseteq M_q $$ and $$W_{q,j}=\mathrm{span}\{\mathbf{x}^{\alpha}:0\leq \alpha_i<q,\ \alpha_1+\dots+\alpha_n = j\} \subseteq M_q $$ for $j \in \mathbb{N}$. For every integer $t \geq 1$, if $ j < \frac{n(q-1)}{2}$, then every k-subspace $W \subseteq V_{q,j}^{\oplus t}$ satisfies $$\dim_k(\sum_{i=1}^n x_iW) \geq \dim_kW.$$ 
\end{lemma}
\begin{proof}
    We first prove the assertion for \(t=1\). First we show the homogeneous case. Assume that $W \subseteq W_{q,j}$ for some $j$. By \cite[Definition 3.2, Proposition 3.4]{dim-ineq}, it suffices to show $\dim_k W_{q,j} \leq \dim_k W_{q,j+1}$, which follows by \cite[Theorem 1, Corollary 2]{dim-formula}.

    For arbitrary $W \subseteq V_{q,j}$, we proceed by induction on \(j\). The case \(j=0\) follows from the homogeneous case. If $j \geq 1$, set $W'=W \cap V_{q,j-1}$. Denote $\pi_d: V_{q,d} \rightarrow W_{q,d}$ the projection onto the homogeneous component of degree \(d\). Then we have $$\sum_{i=1}^n x_iW' \subseteq (\sum_{i=1}^n x_iW) \cap V_{q,j}$$ and $$\pi_{j+1}(\sum_{i=1}^n x_iW)=\sum_{i=1}^n x_i\pi_j(W).$$ Consequently, by the induction hypothesis and the homogeneous case, 
    \begin{equation}
        \begin{aligned}
            \dim_k \sum_{i=1}^n x_iW&=\dim_k ((\sum_{i=1}^n x_iW)\cap V_{q,j})+\dim_k \pi_{j+1}(\sum_{i=1}^n x_iW)\\
            & \geq \dim_k(\sum_{i=1}^n x_iW')+\dim_k(\sum_{i=1}^n x_i\pi_j(W))\\
            & \geq \dim_k W'+\dim_k \pi_j(W)\\
            &=\dim_k W.
        \end{aligned}
    \end{equation}
    The general case follows easily by induction on \(t\).
\end{proof}

\begin{theorem}\label{thm-ha-hb-l-bound}
    Let $(R,\mathfrak{m},k)$ be a local ring of dimension $n$ and characteristic $p>0$ and let $I,J$ be two $\m$-primary ideals, $J^M \subseteq I$, $M \in \mathbb{Z}_{>0}$. Set $$\HH_s(n)=\frac{1}{n!}\sum_{i=0}^{\lfloor s \rfloor}(-1)^i\binom{n}{i}(s-i)^n$$ for $s\geq 0$ and $\HH_s(n)=0$ for $s < 0$, $\HH_s(0)=1$ for $s \geq 0$. Set $$v=\frac{1}{M},\ r=\mu_R(I)+\mu_R((I+J)/I)-n.$$ If $I$ has a reduction $Q$ that is a parameter ideal, we can furthermore choose $r=\mu_R((J+Q)/Q)$. Set $$\Phi_{s,v,r}(n)=\HH_{s}(n)-r\HH_{s-v}(n).$$ Then for $0<a<b\leq v+\frac{n}{2}$, we have
    $$h_b(I,J)-h_a(I,J) \geq e_R(I)(\Phi_{b,v,r}(n)-\Phi_{a,v,r}(n)).$$
\end{theorem}
\begin{proof}
    By Lemma \ref{ring-extension}, we can assume that $R$ is complete and has an infinite residue field. Take $Q=(x_1,\dots,x_n)$ the minimal reduction of $I$, suppose $I^{w+1}=QI^w$. Take $J'=J+Q$. By Proposition \ref{prop-h}\ref{item-h-integral} and Lemma \ref{lemma-hb-ha-J}, we get 
    \begin{equation}\label{h-J-J'-ineq}
        h_b(I,J)-h_a(I,J) \geq h_b(Q,J')-h_a(Q,J').
    \end{equation}
    It suffices to show
    $$h_b(Q,J')-h_a(Q,J') \geq e_R(I)(\Phi_{b,v,r}(n)-\Phi_{a,v,r}(n)). $$

    Take $A=k[[x_1,\dots,x_n]]$, $\n=(x_1,\dots,x_n) \subset A$. Then $\n R=Q$. By Lemma \ref{lemma-rankAR}, writing $\rho=\mathrm{rank}_AR=e_R(Q)=e(I)$, we have an injection $A^\rho \hookrightarrow R$, and there exists $0 \neq c \in A$ such that $cR \subseteq A^\rho$. By Lemma \ref{lemma-J-Q-kq}, there exist $y_1,\dots,y_r \in J$ such that $J'=Q+(y_1,\dots,y_r)$ and $y_i^q \in Q^{k_q}$ for $q$ sufficiently large, where $r=\mu_R((J+Q)/Q),$ $q=p^e,\ k_q=\lfloor q/M \rfloor-w$. Therefore, we have
    \begin{equation}\label{subseteq-chain-A-rho}
        cy_i^qA^\rho \subseteq cQ^{k_q}R = c\n^{k_q}R \subseteq \n^{k_q}A^\rho.
    \end{equation}

    Now we set $$\varphi_q:\ A^{\rho r} \rightarrow A^\rho,\ (w_1,\dots,w_r) \mapsto \sum_{i=1}^r cy_i^qw_i,\ w_i \in A^\rho.$$
    For $t \geq 0,$ we can see that 
    \begin{equation}\label{phi-q-filtration}
        \varphi_q(\n^tA^{\rho r}) \subseteq \n^{t+k_q}A^\rho
    \end{equation}
    by (\ref{subseteq-chain-A-rho}). On the other hand, $A^\rho \hookrightarrow R$ induces
    $$\psi_{q,j}:\ \frac{A^\rho}{(\n^{[q]}+\n^{j+1})A^\rho+\im \varphi_q} \rightarrow \frac{R}{Q^{j+1}+(J')^{[q]}}$$ for each $j \in \mathbb{N}$. This is well-defined since $(J')^{[q]}=Q^{[q]}+(y_1^q,\dots,y_r^q)$ and $\im \varphi_q =\sum cy_i^qA^\rho \subseteq (y_1^q,\dots,y_r^q)R$. Moreover, take $N_q=nq+1,$ then $\n^{N_q} \subseteq \n^{[q]}$, which implies 
    \begin{equation}\label{n^m-M}
        \n^{N_q}\frac{A^\rho}{(\n^{[q]}+\n^{j+1})A^\rho+\im \varphi_q}=0.
    \end{equation}
    We also have
    \begin{equation}\label{n^m-N}
        \n^{N_q} \frac{R}{Q^{j+1}+(J')^{[q]}}=0
    \end{equation}
    since $\n^{N_q}R=Q^{N_q} \subseteq Q^{[q]}$.

    \begin{claim}
        $c \coker \psi_{q,j}=0,\ c^2\ker \psi_{q,j}=0$. 
    \end{claim}
    \begin{proof}
        $c \coker \psi_{q,j}=0$ follows from $cR \subseteq A^\rho$. Let $z \in \ker \psi_{q,j}$. Then $z \in Q^{j+1}+(J')^{[q]}=(\n^{j+1}+\n^{[q]})R+(y_1^q,\dots,y_r^q)R.$ Write $z=\sum_v a_vr_v+\sum_{i=1}^r y_i^qr_i,\ a_v \in \n^{j+1}+\n^{[q]},\ r_v,r_i \in R$. Since $cr_i,cr_v \in A^\rho,$ we have 
        \begin{equation}
            \begin{aligned}
                c^2z &=\sum_v a_v(c^2r_v)+\sum_{i=1}^r cy_i^q(cr_i)\\
                &=\sum_v a_v(c^2r_v) +\varphi_q(cr_1,\dots,cr_r) \\
                &\in (\n^{j+1}+\n^{[q]})A^\rho+\im \varphi_q.
            \end{aligned}
        \end{equation}
        Therefore, $c^2\ker \psi_{q,j}=0$.  
    \end{proof}
    Note that $$\mu_A(\frac{R}{Q^{j+1}+(J')^{[q]}}) \leq \mu_A(R),\ \mu_A(\frac{A^\rho}{(\n^{[q]}+\n^{j+1})A^\rho+\im \varphi_q}) \leq \rho.$$    
    Combining (\ref{n^m-M}) and (\ref{n^m-N}), by Lemma \ref{lemma-length-eq} and Lemma \ref{lemma-rankAR}, we know that 
    \begin{equation}\label{length-A-j=length-R-j}
        \lambda_A(\frac{A^\rho}{(\n^{[q]}+\n^{j+1})A^\rho+\im \varphi_q})=\lambda_R(\frac{R}{Q^{j+1}+(J')^{[q]}})+O(q^{n-1}).
    \end{equation}
    In fact, if \(c\) is a unit, then \(\psi_{q,j}\) is an isomorphism; otherwise, we apply Lemma \ref{lemma-length-eq} with \(c^2\).

    By the definition of $\varphi_q$, we have $\varphi_q(\n^{[q]}A^{r\rho}) \subseteq \n^{[q]}A^\rho,$ which induces
    $$\overline{\varphi_q}:\ (A/\n^{[q]})^{\rho r} \rightarrow (A/\n^{[q]})^\rho. $$ We denote $$M_q=A/\n^{[q]}=\frac{k[[x_1,\dots,x_n]]}{(x_1^q,\dots,x_n^q)}.$$ A $k$-basis of $M_q$ is given by monomials $\mathbf{x}^{\alpha}$, where $\alpha=(\alpha_1,\dots,\alpha_n) \in \mathbb{N}^n, 0 \leq \alpha_i <q.$ Denote $$V_{q,j}=\mathrm{span}\{\mathbf{x}^{\alpha}:0\leq \alpha_i<q,\ \alpha_1+\dots+\alpha_n \leq j\}=\frac{A}{\n^{j+1}+\n^{[q]}} \subseteq M_q$$ and $V_{q,j}=0$ for $j<0$. Denote the natural projection $\pi_{q,j}:M_q \twoheadrightarrow V_{q,j}$. By (\ref{phi-q-filtration}), we denote $\varphi_{q,j}: V_{q,j-k_q}^{\rho r} \rightarrow V_{q,j}^{\rho}$ the morphism induced by $\overline{\varphi_q}$.
    \begin{claim}
        $\pi_{q,j}\circ\overline{\varphi_q}((A/\n^{[q]})^{\rho r}) = \im \varphi_{q,j}.$
    \end{claim}
    \begin{proof}
        If $j<k_q$, then $\varphi_q(A^{\rho r}) \subseteq \n^{k_q}A^\rho \subseteq \n^{j+1}A^\rho.$ Thus the image after projection is zero. If $j\geq k_q$, the result follows from the commutative diagram of $k$-linear maps

$$\begin{tikzcd}
M_q^{\rho r} \arrow[r, "\overline{\varphi_q}"] \arrow[d, "{\pi_{q,j-k_q}}"] & M_q^{\rho} \arrow[d, "{\pi_{q,j}}"] \\
{V_{q,j-k_q}^{\rho r}} \arrow[r, "{\varphi_{q,j}}"]                         & {V_{q,j}^{\rho}}.                  
\end{tikzcd}$$ 
    \end{proof}
    This induces $$\frac{A^\rho}{(\n^{[q]}+\n^{j+1})A^\rho+\im \varphi_q} \cong \frac{M_q^\rho}{\im \overline{\varphi_q}+\n^{j+1}M_q^\rho} \cong \frac{V_{q,j}^\rho}{\im \varphi_{q,j}}.$$
    Combining (\ref{length-A-j=length-R-j}), we have
    \begin{equation}\label{length-dim-eq}
        \lambda_R(\frac{R}{Q^{j+1}+(J')^{[q]}})=\rho \dim_kV_{q,j}-\rho r \dim_k V_{q,j-k_q}+\dim_k \ker \varphi_{q,j} +O(q^{n-1}).
    \end{equation}

    For the dimension of the kernel, if $j-k_q<0$, then $\ker \varphi_{q,j}=0$. Now assume $0 \leq j-k_q<\frac{n(q-1)}{2}$. For every $t \in \ker \varphi_{q,j}$, 
    \begin{equation}
        \begin{aligned}
            \varphi_{q,j+1}(x_it)&= \pi_{q,j+1}\circ \overline{\varphi_{q}}(x_it)\\
            &=\pi_{q,j+1}(x_i \overline{\varphi_q}(t))\\
            &=x_i\pi_{q,j}(\overline{\varphi_q}(t))\\
            &=x_i\varphi_{q,j}(t)=0.
        \end{aligned}
    \end{equation}
    Therefore, $\sum_{i=1}^n x_i\ker\varphi_{q,j} \subseteq \ker \varphi_{q,j+1}$. Applying Lemma \ref{lemma-dim-ineq} on $\ker \varphi_{q,j} \subseteq V_{q,j-k_q}^{\rho r}$, we get
    \begin{equation}\label{ker-dim-ineq}
        \dim_k \ker \varphi_{q,j+1} \geq \dim_k \ker\varphi_{q,j}.
    \end{equation}

    For the dimension of $V_{q,j}$, fix $u>0$, we have $$\lim_{e \rightarrow \infty} \frac{1}{q^n}\dim_k V_{q,\lceil uq \rceil-1}=\lim_{e \rightarrow \infty} \frac{1}{q^n} \lambda_A(\frac{A}{\n^{\lceil uq \rceil}+\n^{[q]}})=\HH_u(n)$$ by \cite[Proposition 3.4]{e-and-eHK}, where $e \rightarrow \infty$ means that $q=p^e \rightarrow \infty$. It easily follows that for any integer sequence $u(q)$ such that $u(q)/q \rightarrow u$ as $e \rightarrow \infty$,
    \begin{equation}\label{dimV-H}
        \lim_{e \rightarrow \infty} \frac{1}{q^n}\dim_k V_{q,u(q)}=\HH_u(n).
    \end{equation}

    Now for $0 < a<b<v+\frac{n}{2}$, denote $j_a=\lceil aq \rceil-1,\ j_b=\lceil bq \rceil-1.$ Recall that $k_q=\lfloor \frac{q}{M} \rfloor-w$ and $v=\frac{1}{M}$, we get $\frac{j_b-k_q}{q} \rightarrow b-v<\frac{n}{2}$ as $e \rightarrow \infty$. Therefore for sufficiently large $q$, $j_b-k_q<\frac{n(q-1)}{2}$. Applying (\ref{ker-dim-ineq}), we get 
    $$\dim_k \ker \varphi_{q,j_b} \geq \dim_k \ker\varphi_{q,j_a}.$$ Applying (\ref{length-dim-eq}) on $j_a,j_b$, we get
    \begin{equation}
        \begin{aligned}
            \lambda_R(\frac R{Q^{\lceil bq\rceil}+(J')^{[q]}})-\lambda_R(\frac R{Q^{\lceil aq\rceil}+(J')^{[q]}}) \geq &\rho(\dim_kV_{q,j_b}-\dim_kV_{q,j_a})\\&-\rho r(\dim_kV_{q,j_b-k_q}-\dim_kV_{q,j_a-k_q})+O(q^{n-1}).
        \end{aligned}
    \end{equation}
    
    Dividing both sides by $q^n$, combining (\ref{dimV-H}), (\ref{h-J-J'-ineq}) and $\rho=e_R(I)$, we find that
    \begin{equation}\label{main-ineq-half}
        h_b(I,J)-h_a(I,J) \geq e_R(I)(\Phi_{b,v,r}(n)-\Phi_{a,v,r}(n))
    \end{equation}
    holds for $r=\mu_R((J+Q)/Q)$ and $0 < a<b<v+\frac{n}{2}$. Applying (\ref{main-ineq-half}) to $b_l \rightarrow b,$ we can show (\ref{main-ineq-half}) holds for $b=v+\frac{n}{2}$. Applying (\ref{main-ineq-half}) to ideals $I,I+J$ and combining Lemma \ref{lemma-hb-ha-J}, we have $h_b(I,J)-h_a(I,J) \geq h_b(I,I+J)-h_a(I,I+J) \geq e_R(I)(\Phi_{b,v,r_0}(n)-\Phi_{a,v,r_0}(n))$ for $r_0=\mu_R((I+J)/Q)$. We have $$0 \rightarrow I/Q \rightarrow (I+J)/Q \rightarrow (I+J)/I \rightarrow 0, $$ which implies $$r_0=\mu_R((I+J)/Q) \leq \mu_R(I/Q)+\mu_R((I+J)/I)=\mu_R(I)+\mu_R((I+J)/I)-n.$$ Therefore, (\ref{main-ineq-half}) holds for $r=\mu_R(I)+\mu_R((I+J)/I)-n$. This finishes the proof.
\end{proof}

To extract a quadratic lower bound from Theorem \ref{thm-ha-hb-l-bound}, we study the derivative of \(\Phi_{s,v,r}(n)\). We first recall several elementary properties of \(\HH_s(n)\).

\begin{lemma}\label{lemma-H_s}
    Let $n\geq 0$, $s$ be a real number. Let $\HH_n(s)$ be defined as in Theorem \ref{thm-ha-hb-l-bound}. Then we have
    \begin{enumerate}
        \item\label{item-H-volume} $\HH_s(n)$ is the volume of the set $$\{(x_1,\dots,x_n) \in [0,1]^n : x_1+\dots+x_n \leq s\} \subset \mathbb{R}^n.$$
        \item\label{item-H-scaling}  For $s \geq 1,\ n \geq 1$, we have $$\frac{s^{n}}{n!}-n\frac{(s-1)^{n}}{n!} \leq \HH_s(n)\leq \frac{s^{n}}{n!} .$$
        \item\label{item-H-diff} For $n \geq 1$, $\mathcal{H}_s(n)=\int_{s-1}^{s}\mathcal{H}_{t}(n-1)dt$. Therefore, for $n \geq 2,$ $\HH_s'(n)=\HH_s(n-1)-\HH_{s-1}(n-1).$
        \item\label{item-H'-scaling} For $s \geq 1,\ n \geq 2$, we have $$\frac{s^{n-1}}{(n-1)!}-n\frac{(s-1)^{n-1}}{(n-1)!} \leq \HH_s’(n)\leq \frac{s^{n-1}}{(n-1)!} .$$
    \end{enumerate}
\end{lemma}
\begin{proof}
    (i) For each integer $0 \leq i \leq \min\{\lfloor s \rfloor,n\} $, $\frac{1}{n!}(s-i)^n$ is exactly the volume of the set $$\{(x_1,\dots,x_n) \in \mathbb{R}_{\geq 0}^{n}: x_1+\dots+x_n \leq s,\ x_1,\dots,x_i \geq 1\}. $$ By the inclusion–exclusion principle, $\HH_s(n)$ equals the volume of $$\{(x_1,\dots,x_n) \in [0,1]^n : x_1+\dots+x_n \leq s\} .$$\\
    (ii) We have $$\frac{s^{n}}{n!}-n\frac{(s-1)^{n}}{n!} \leq \HH_s(n)\leq \frac{s^{n}}{n!} $$ by \ref{item-H-volume} and the inclusion–exclusion principle.\\
    (iii) By \cite[Lemma 3.7]{e-and-eHK}.\\
    (iv) Let $t>0$. By \ref{item-H-volume}, $\HH_{s+t}(n)-\HH_s(n)$ is the volume of the set $$\{(x_1,\dots,x_n) \in [0,1]^n : s \leq x_1+\dots+x_n \leq s+t\} .$$ For each integer $0 \leq i \leq\min\{\lfloor s \rfloor,n\}$, $\frac{(s+t-i)^n-(s-i)^n}{n!}$ is the volume of the set $$\{(x_1,\dots,x_n) \in \mathbb{R}_{\geq 0}^n : s \leq x_1+\dots+x_n \leq s+t,\ x_1,\dots,x_i \geq 1 \}.$$ Therefore, by the inclusion–exclusion principle, we have 
    $$\frac{(s+t)^n-s^n-n((s+t-1)^n-(s-1)^n)}{n!} \leq \HH_{s+t}(n)-\HH_s(n) \leq \frac{(s+t)^n-s^n}{n!}.$$ Dividing by $t$ and letting $t \rightarrow 0^+$, we get $$\frac{s^{n-1}}{(n-1)!}-n\frac{(s-1)^{n-1}}{(n-1)!} \leq \HH_s’(n)\leq \frac{s^{n-1}}{(n-1)!}$$ for $s \geq 1,\ n \geq 2$.
\end{proof}

\begin{lemma}\label{lemma-tau-exist}
    Let $r \geq 1,  n \geq 2$ be integers, $0<v<n$ be a real number, $\HH_s(n),\ \Phi_{s,v,r}$ be defined as in Theorem \ref{thm-ha-hb-l-bound}, $s \in \mathbb{R}$. Then $\Phi_{s,v,r}'(n)$ has a unique zero $\tau$ in the interval $(v,n)$, where $\Phi_{s,v,r}'(n)$ denotes the derivative of \(\Phi_{s,v,r}(n)\) with respect to \(s\), and $v < \tau \leq \frac{v+n}{2}$. If $r=1$, then $\tau=\frac{v+n}{2}$. Moreover, $\tau$ is a maximum point of $\Phi_{s,v,r}(n)$ and $$\Phi_{\tau,v,r}''(n)<0 \quad \mathrm{if\ }n \geq 3.$$ 
\end{lemma}
\begin{proof}
    The proof is similar to \cite[Lemma 3.9]{lower-bound-s}. Write \(F_n(s)=\Phi'_{s,v,r}(n)\). We prove by induction on \(n\ge2\) to show that \(F_n\) changes sign exactly once, from positive to negative, on \((0,n+v)\). The case \(n=2\) follows directly from $\mathcal H'_s(2)=\max\{0,\min\{s,2-s\}\}.$

Let \(n\ge3\). Since \(F_n>0\) on \((0,v]\) and \(F_n<0\) on \([n,n+v)\), it has a zero in \((v,n)\). If \(n-1\le v<n\), then, for \(v<s<n\),
\[
F_n(s)=\frac{(n-s)^{n-1}-r(s-v)^{n-1}}{(n-1)!},
\]
whose derivative is strictly negative, proving all assertions.

Suppose \(v<n-1\), and let \(\sigma\) be the unique zero of \(F_{n-1}\) in \((v,n-1)\). By Lemma \ref{lemma-H_s},
\[
F_n(s)=\int_{s-1}^{s}F_{n-1}(t)\,dt,
\qquad
F'_n(s)=F_{n-1}(s)-F_{n-1}(s-1).
\]
The sign change of \(F_{n-1}\) implies that every zero \(\tau\in(v,n)\) of \(F_n\) satisfies $\tau-1<\sigma<\tau.$ Moreover, for \(\sigma<s<\sigma+1\), we have
\[
F_{n-1}(s)\le0\le F_{n-1}(s-1),
\]
with at least one strict inequality: otherwise \(s\ge n-1+v\) and \(s\le1\), a contradiction. Thus \(F'_n(s)<0\) throughout \((\sigma,\sigma+1)\). Consequently, the zero \(\tau\) is unique, $\Phi''_{\tau,v,r}(n)=F'_n(\tau)<0,$ and the sign change of \(F_n\) makes \(\tau\) the unique global maximizer of \(\Phi_{s,v,r}(n)\).

By symmetry, \(\mathcal H'_s(n)=\mathcal H'_{n-s}(n)\), and hence
\[
\Phi'_{\frac{n+v}{2},v,r}(n)
=(1-r)\mathcal H'_{\frac{n-v}{2}}(n)\le0.
\]
The sign change of \(\Phi'_{s,v,r}(n)\) therefore gives
\(\tau\le(n+v)/2\), with equality when \(r=1\).

\end{proof}

\begin{corollary}\label{cor-hb-ha-qua}
    Retain the notation of Theorem \ref{thm-ha-hb-l-bound} and suppose $n \geq 2$, $r \geq 1$. Then $\Phi_{s,v,r}'(n)=0$ has a unique solution $\tau$ in the interval $(v,n)$ and $0<v<\tau\leq\frac{n+v}{2}<v+\frac{n}{2}$. For any fixed real number \(0<\delta<\tau\), there exists a constant $C>0$ such that $$h_{b}(I,J)-h_a(I,J) \geq C(b-a)^2$$ for every $\delta \leq a \leq b \leq \tau$. 
    
    In particular, if $I=J$ and $\mu_R(I) \leq n+1$, we may take $v=r=1$. Then $\tau=\frac{n+1}{2}$ and $h_{b}(I,I)-h_a(I,I) \geq C(b-a)^2$ for $\delta \leq a \leq b \leq \frac{n+1}{2}$.
\end{corollary}
\begin{proof}
    The existence and the stated bounds of \(\tau\) are proved in Lemma \ref{lemma-tau-exist}. In particular, if $I=J$ and $\mu_R(I) \leq n+1$, we may take $v=1,\ r=1$. Lemma \ref{lemma-tau-exist} gives $\tau=\frac{n+1}{2}$. 

    Lemma \ref{lemma-tau-exist} also shows that $\Phi_{\tau,v,r}''(n)<0 $ for $n \geq 3$. For $n=2$, an easy calculation shows that the left-hand limit of $\Phi_{s,v,r}''(n)$ at $\tau$ is negative. Therefore, for every $0<\delta<\tau$, there exists $c>0$ such that $\Phi_{s,v,r}'(n)>c(\tau-s)$ for $\delta\leq s<\tau$.

    For $\delta\leq a \leq b \leq \tau$, by Theorem \ref{thm-ha-hb-l-bound}, we have
    \begin{equation}
        \begin{aligned}
            h_b(I,J)-h_a(I,J) &\geq e_R(I) (\Phi_{b,v,r}(n)-\Phi_{a,v,r}(n))\\
            &=e_R(I)\int_a^b \Phi_{s,v,r}'(n) ds\\
            &\geq e_R(I) \int_a^b c(\tau-s)ds\\
            &=e_R(I)c(\tau(b-a)-\frac{b^2-a^2}{2})\\
            &=e_R(I)c((b-a)(\tau-b)+\frac{(b-a)^2}{2})\\
            &\geq \frac{ce_R(I)}{2}(b-a)^2.
        \end{aligned}
    \end{equation}
\end{proof}

\section{Applications}\label{section-application}
\subsection{Lower bounds for $F$-thresholds}
As an application, for a Noetherian local ring of positive characteristic with dimension not less than $3$, we can give a lower bound of $F$-thresholds.

\begin{theorem}\label{thm-cJI}
    Let $(R,\m)$ be a Noetherian local ring of dimension $n \geq 3$ and characteristic $p>0$, and let $I,J$ be two $\m$-primary ideals, $J^M \subseteq I$, $M \in \mathbb{Z}_{>0}$. Set $$v=\frac{1}{M},\ r=\mu_R(I)+\mu_R((I+J)/I)-n,\ \Phi_{s,v,r}(n)=\HH_{s}(n)-r\HH_{s-v}(n).$$  If $I$ has a reduction $Q$ that is a parameter ideal, we can furthermore choose $r=\mu_R((J+Q)/Q)$. If $r \geq 1$, let $\tau$ be the unique solution in $(v,n)$ of $\Phi_{\tau,v,r}'(n)=0$. Then $c^J(I) \geq \tau$. Moreover, if 
    \begin{equation}\label{condition-mu-omega}
        \sup_{q=p^e} \mu_R(\omega_{R/J^{[q]}})<\infty,
    \end{equation}
    then $c^J(I)>\tau.$

    In particular, we have
    \begin{equation}\label{cJI-ineq-eq}
        c^J(I) \geq \begin{cases}
            \dfrac{n+1}{2M}, &\mu_R(J) \leq n+1,\\
            \dfrac{\mu_R(J)^{1/(n-1)}}{M(\mu_R(J)^{1/(n-1)}-1)},&\mu_R(J)\ge n+2.
        \end{cases}
    \end{equation}
    If (\ref{condition-mu-omega}) is satisfied, then the sign $\geq$ in (\ref{cJI-ineq-eq}) can be replaced by $>$.
\end{theorem}
\begin{proof}
    Assume that $c^J(I) < \tau.$ Choose $\epsilon>0$ sufficiently small such that $c^J(I)<\tau-\epsilon$. By the definition of $c^J(I)$, we have $I^{\lceil (\tau-\epsilon) q \rceil} \subseteq J^{[q]}$ for sufficiently large $q=p^e$. Thus $h_{\tau-\epsilon}(I,J)=h_{\tau}(I,J)=e_{HK}(J)$, which leads to a contradiction to Corollary \ref{cor-hb-ha-qua}.

    Now we assume that (\ref{condition-mu-omega}) is satisfied and set $B=\sup_{q=p^e} \mu_R(\omega_{R/J^{[q]}})$. If $c^J(I) \leq \tau,$ the first part implies $c^J(I)=\tau$. We take $N_q=\nu^J_I(q)+1$. Then $I^{N_q} \subseteq J^{[q]}$ for sufficiently large $q=p^e$. Fix $\epsilon>0$ sufficiently small, denote $a_q=\lfloor(\tau-\epsilon)q\rfloor,\ b_q=N_q-a_q$. Then $a_q,b_q>0,\ I^{a_q+b_q} \subseteq J^{[q]}$ and $a_q/q \rightarrow \tau-\epsilon,\ b_q/q \rightarrow \epsilon$ as $q=p^e \rightarrow \infty$.

    Applying Lemma \ref{lemma-long-ineq} to $I,J^{[q]}$, we have
    \begin{equation}
        \lambda_R(\frac{I^{a_q}+J^{[q]}}{I^{N_q}+J^{[q]}}) = \lambda_R(\frac{I^{a_q}+J^{[q]}}{J^{[q]}}) \leq \mu_R(\omega_{R/J^{[q]}})\lambda_R(\frac{R}{I^{b_q}+J^{[q]}}) \leq B \lambda_R(\frac{R}{I^{b_q}}).
    \end{equation}
    Dividing both sides by $q^n$ and taking the limit, we get
    \begin{equation}
        h_\tau(I,J)-h_{\tau-\epsilon}(I,J) \leq C'\epsilon^n,
    \end{equation}
    where $C'=\frac{Be(I)}{n!}>0$ is a constant independent of $\epsilon$. On the other hand, we have $h_\tau(I,J)-h_{\tau-\epsilon}(I,J) \geq C \epsilon^2$ for some $C>0$ by Corollary \ref{cor-hb-ha-qua}. Since $n \geq 3$, there is a contradiction, which shows that $c^J(I)>\tau$.

    We next give an estimate of $\tau$. We only need to consider $\tau$ for the ideal pair $(J,J)$, hence $r=\mu_R(J)-n,\ v=1$. Then $c^J(I) \geq c^J(J^M)=\frac{1}{M}c^J(J)$. Clearly $\mu_R(J) \geq n$. If $\mu_R(J)=n$, $J$ is a parameter ideal and $c^J(J)=n>\frac{n+1}{2}$ by \cite[Example 2.7(i)]{HMTW}. If $\mu_R(J)=n+1$, $\tau=\frac{n+1}{2}$ by \cite[Lemma 3.9]{lower-bound-s}.

    For $\mu_R(J) \geq n+2$, we consider $\Phi_{s,v,r}'(n)=\HH_s'(n)-r\HH_{s-v}'(n)=\HH_s'(n)-(\mu_R(J)-n)\HH_{s-1}'(n)$. By Lemma \ref{lemma-H_s}\ref{item-H'-scaling}, we have 
    \begin{equation}\label{Phi'-l-bound}
        \Phi_{s,v,r}'(n) \geq \frac{s^{n-1}-\mu_R(J)(s-1)^{n-1}}{(n-1)!}
    \end{equation}
    for $s \geq 1$. Take $$\tau_0=\frac{\mu_R(J)^{1/(n-1)}}{\mu_R(J)^{1/(n-1)}-1}>1,$$ then $\Phi_{s,v,r}'(n) \geq 0$ for $0<s<\tau_0$ ($1 \leq s <\tau_0$ follows from (\ref{Phi'-l-bound}), $0 < s \leq 1$ is trivial). Therefore, $\tau \geq \tau_0$, and $c^J(J) \geq \tau \geq \tau_0$. It follows that $c^J(I) \geq \frac{1}{M}\tau_0$.
\end{proof}
\begin{remark}
    The bound $\tau$ is sharp, and in general $\geq$ cannot be replaced by $>$. See the following example.

    Take $R=k[[x,y,z,w]]/(xy-zw)$ with $\chara k=p>0$, $I=J=\m=(x,y,z,w)$. Then $n=\dim R=3$ and $\mu_R(I)=4$. For every $q=p^e$, one can show $\m^{2q-1} \subseteq \m^{[q]},\ \m^{2q-2} \nsubseteq \m^{[q]}$. Therefore, $c^\m(\m)=2=\frac{n+1}{2}$.
\end{remark}

To apply the strict part of Theorem \ref{thm-cJI}, we verify the required uniform bound in two familiar settings.

\begin{proposition}\label{prop-condition-hold}
    Condition (\ref{condition-mu-omega}) holds in each of the following cases. 
    \begin{enumerate}
        \item\label{item-condition-regular} $R$ is regular.
        \item $R$ is Cohen-Macaulay and $J$ is a parameter ideal.
    \end{enumerate}
    In fact, in either case, $\mu_R(\omega_{R/J^{[q]}})=\mu_R(\omega_{R/J})$ for every $q=p^e$.
\end{proposition}
\begin{proof}
    (i) If $R$ is regular, the Frobenius map is flat. Denote $r(M)=\dim_k(0:_M \m)$ the same as in Lemma \ref{lemma-duality}\ref{item-canonical}. By Lemma \ref{lemma-duality}\ref{item-canonical} and \cite[Proposition 1.2.16(b)]{cohen-macaulay-ring}, we have $$\mu_R(\omega_{R/J^{[q]}})=r(R/J^{[q]})=r(R/J)r(R/\m^{[q]})=r(R/J)=\mu_R(\omega_{R/J}).$$ $R/\m^{[q]}$ is an Artinian complete intersection, hence Artinian Gorenstein. Lemma \ref{lemma-duality}\ref{item-gorenstein} gives $r(R/\m^{[q]})=1$.\\
    (ii) If $R$ is Cohen-Macaulay and $J=(x_1,\dots,x_n)$ is a parameter ideal, then $x_1^q,\dots,x_n^q$ is a regular sequence. Therefore $\mu_R(\omega_{R/J^{[q]}})=r(R/J^{[q]})=r(R/J)=\mu_R(\omega_{R/J})$ by \cite[Lemma 1.2.19]{cohen-macaulay-ring}.
\end{proof}

As a corollary, we directly get
\begin{corollary}\label{cor-cII}
    Let $(R,\m)$ be a regular Noetherian local ring of dimension $n \geq 3$ and characteristic $p>0$, and let $I$ be an $\m$-primary ideal, $\mu_R(I)\leq n+1$. Then $c^I(I)>\frac{n+1}{2}$. 
\end{corollary}

\subsection{A result about $s$-closure}

In \cite{e-and-eHK}, William D. Taylor introduced the definition of $s$-closures, which are closely related to $h$-functions. In this part, we will prove Theorem \ref{intro-thm-nsubset-s-closure}, a noncontainment statement for ideal quotients and \(s\)-closures. 

\begin{definition}\textup{(\cite[Definition 4.1, Definition 4.3]{e-and-eHK})}
    Let $R$ be a ring of characteristic $p>0$, let $I$ be an ideal of $R$, and let $s \geq 1$ be a real number. 
    
    An element $x \in R$ is said to be in the weak $s$-closure of $I$ if there exists $c \in R^\circ$ such that for all $q=p^e$ sufficiently large, $cx^q \in I^{\lceil sq \rceil}+I^{[q]},$ where $R^\circ$ stands for the complement of the union of the minimal primes of $R$. We denote the set of all $x$ in the weak $s$-closure of $I$ by $I^{\wcl_s}$.

    The $s$-closure of $I$ is defined to be the union of the following chain of ideals:
    $$I \subseteq I^{\text{w.cl}_s} \subseteq (I^{\text{w.cl}_s})^{\text{w.cl}_s} \subseteq ((I^{\text{w.cl}_s})^{\text{w.cl}_s})^{\text{w.cl}_s} \subseteq \dots .$$ We denote this ideal by $I^{\text{cl}_s}$.
\end{definition}

First we give a formula on the quotient of Frobenius power and general power of parameter ideals.

\begin{lemma}\label{lemma-Qq:Qa}
    Let $R=k[[x_1,\dots,x_n]]$, where $\chara k=p>0$. Let $Q$ be a parameter ideal of $R$. Then for every $q=p^e$ and integer $a \geq 0$, we have $$Q^{[q]}:Q^a=Q^{[q]}+Q^{n(q-1)+1-a},$$ where $Q^b$ is regarded as $R$ if $b<0$.
\end{lemma}
\begin{proof}
    Write $Q=(y_1,\dots,y_n)$. Then $y_1,\dots,y_n$ is a system of parameters. Let $A=k[[y_1,\dots,y_n]]$, $\n=(y_1,\dots,y_n)$ be the maximal ideal of $A$. By Lemma \ref{lemma-rankAR}, $R$ is a finite $A$-module, and $\n^{[q]}R=Q^{[q]},\ \n^a R=Q^a$. By \cite[Proposition 7.4]{tight-closure}, we have
    \begin{equation}\label{Q-n-quotient-monomial}
        Q^{[q]}:Q^{a}=(\n^{[q]}:\n^a)R.
    \end{equation}

    $\n^{[q]},\n^a$ are monomial ideals, thus we only need to consider the monomial $\mathbf{y}^\alpha$ for $\alpha=(\alpha_1,\dots,\alpha_n)$. $\mathbf{y}^\alpha \notin \n^{[q]}$ if and only if $0 \leq \alpha_i \leq q-1.$ Denote $D=n(q-1)$, $|\alpha|=\alpha_1+\dots+\alpha_n$.

    If $|\alpha| \geq D+1-a$, for any $|\beta|=a$, $|\alpha+\beta| \geq D+1$, there exists an $i$ such that $\alpha_i+\beta_i \geq q$, which implies $\mathbf{y}^{\alpha+\beta} \in \n^{[q]}$. That is, $\mathbf{y}^\alpha \in \n^{[q]}:\n^a$.

    If $|\alpha| \leq D-a$, take $\mathbf{y}^{\gamma}=\Pi_{i=1}^n y_i^{q-1-\alpha_i}$. Then $|\gamma| =D-|\alpha| \geq a$ and $\mathbf{y}^{\gamma} \in \n^a$. But $\mathbf{y}^{\alpha+\gamma}=\Pi_{i=1}^n y_i^{q-1} \notin \n^{[q]}$, thus $\mathbf{y}^\alpha \notin \n^{[q]}:\n^a$.

    Therefore, $\n^{[q]}:\n^a=\n^{[q]}+\n^{D+1-a}$. Returning to (\ref{Q-n-quotient-monomial}), we have $Q^{[q]}:Q^a=Q^{[q]}+Q^{n(q-1)+1-a} $.
\end{proof}

Then we give a characterization of \(s\)-closure in formal power series rings.

\begin{proposition}\label{prop-cl-wcl-T}
    Let $R=k[[x_1,\dots,x_n]]$, where $\chara k=p>0$. Let $Q$ be a parameter ideal of $R$. For $1 \leq s \leq n$, we have 
    $$Q^{\cl_s}=Q^{\wcl_s}=\{x \in R:x^qQ^{\lceil (n-s)q \rceil} \subseteq  Q^{[q]} \mathrm{\ for\ all\ }q=p^e,\ e \geq 0\}.  $$
\end{proposition}
\begin{proof}
    Denote $T_s=\{x \in R:x^{q}Q^{\lceil (n-s)q \rceil} \subseteq  Q^{[q]} \mathrm{\ for\ all\ }q=p^e,\ e \geq 0\}.$ First we show $Q^{\wcl_s} \subseteq T_s$.

    Take $x \in Q^{\wcl_s}$, for any $z \in Q^{\lceil (n-s)q \rceil}$ and sufficiently large $p$-power $P$, we have
    \begin{equation}
        \begin{aligned}
            c(x^qz)^P &\in (Q^{\lceil sqP \rceil}+Q^{[qP]}) \cdot Q^{\lceil (n-s)q \rceil P} \\&\subseteq Q^{\lceil sqP \rceil +\lceil (n-s)q \rceil P} + Q^{[qP]}\\
            & \subseteq Q^{nqP}+Q^{[qP]}=Q^{[qP]}=(Q^{[q]})^{[P]}
        \end{aligned}
    \end{equation}
    for some $0 \neq c \in R$. This implies $x^qz \in (Q^{[q]})^*,$ the tight closure of $Q^{[q]}$, which coincides with $Q^{[q]}$ in a regular local ring. Therefore, $x \in T_s$.

    Conversely, take $x \in T_s$. By Lemma \ref{lemma-Qq:Qa}, $x^q \in Q^{[q]}+Q^{n(q-1)+1-\lceil (n-s)q \rceil}.$ Choosing $0 \neq c_0 \in Q^n,$ we have $$c_0x^q \in Q^{[q]}+Q^{nq+1-\lceil (n-s)q \rceil}=Q^{[q]}+Q^{\lfloor sq \rfloor+1} \subseteq Q^{[q]}+Q^{\lceil sq \rceil}.$$ This shows $Q^{\wcl_s}=T_s$. It suffices to show that  $Q^{\cl_s} \subseteq Q^{\wcl_s}$.

    If $s=1$, $Q^{\wcl_s}$ and $Q^{\cl_s}$ coincide with the integral closure of $Q$. If $s>1$, take $x \in Q^{\cl_s}$. Fix $q=p^e$ and set $u_P=s-\frac{1}{qP}$ for sufficiently large $p$-power $P$, then $1 \leq u_P<s.$ By \cite[Theorem 4.11]{result-s-closure}, $x \in Q^{\wcl_{u_P}}=T_{u_P}$. Thus, $$x^{qP}Q^{\lceil (n-u_P)qP \rceil}=x^{qP}Q^{\lceil (n-s)qP \rceil +1} \subseteq Q^{[qP]}.$$ On the other hand, taking $z \in Q^{\lceil (n-s)q \rceil}$ and $0 \neq c \in Q$, we have $$cz^P\in Q^{1+P \lceil (n-s)q \rceil }\subseteq Q^{\lceil (n-s)qP \rceil +1}.$$ Therefore, for every $P$ sufficiently large, $c(x^qz)^P \in Q^{[qP]}=(Q^{[q]})^{[P]}$. Since $Q^{[q]}$ is tightly closed, $x^qz \in Q^{[q]}$, which implies $Q^{\cl_s} \subseteq T_s=Q^{\wcl_s}$.   
\end{proof}

\begin{theorem}\label{thm-nsubset-s-closure}
    Let $(R,\m)$ be a complete regular local ring of dimension $n \geq 3$ and characteristic $p>0$. Let $I$ be an $\m$-primary ideal of $R$, with a reduction $Q$ that is a parameter ideal. Let $r=\mu_R(I)-n \geq 1.$ Let $1<\tau \leq \frac{n+1}{2}$ be the unique solution of $\Phi'_{s,1,r}(n)=\HH_s'(n)-r\HH_{s-1}'(n)=0$ in $(1,n)$ as defined in Corollary \ref{cor-hb-ha-qua}. Then $Q:I \nsubseteq Q^{\cl_{n-\tau}}$.
\end{theorem}

\begin{proof}
    Assume that $ Q:I \subseteq Q^{\cl_{n-\tau}}$. By Proposition \ref{prop-cl-wcl-T}, we have $(Q:I)^{[q]} Q^{\lceil \tau q \rceil}\subseteq Q^{[q]}$ for $q=p^e$. Therefore, 
    \begin{equation}\label{J-J+f}
        Q^{\lceil \tau q \rceil} \subseteq Q^{[q]}:(Q:I)^{[q]}=(Q:(Q:I))^{[q]}=I^{[q]}
    \end{equation}
    by flatness of Frobenius map and Lemma \ref{lemma-duality}\ref{item-0:0:I}. (\ref{J-J+f}) implies $c^I(Q)\leq \tau$. Since $\overline{I}=\overline{Q}$, $c^I(I)=c^I(Q) \leq \tau$ by \cite[Proposition 2.2(2)]{HMTW}, which contradicts Theorem \ref{thm-cJI} and Proposition \ref{prop-condition-hold}. That is, $Q:I \nsubseteq Q^{\cl_{n-\tau}}$.
\end{proof}

Therefore, we have answered \cite[Conjecture 3.12]{mu-tau} by applying Theorem \ref{thm-nsubset-s-closure} to the Tjurina ideal and the Jacobian ideal of an isolated hypersurface singularity $f$. That is,

\begin{corollary}\label{cor-mu-s-closure}
    Let $f \in k[[x_1,\dots,x_n]]$ define an isolated hypersurface singularity, $n \geq 3$, $\chara k>0$. Let $j(f)$ denote its Jacobian ideal. Assume that $f\in \overline{j(f)}$. Then $j(f):f \nsubseteq j(f)^{\mathrm{cl}_{\frac{n-1}{2}}}$.
\end{corollary}
\begin{proof}
    If $f \in j(f)$, the result holds automatically. For $f \notin j(f)$, take $I=(f,j(f))$, $Q=j(f)$. If $\mu(I)=n,$ then $I$ is a parameter ideal and hence $\lambda_R(R/Q)=e(Q)=e(I)=\lambda_R(R/I),$ which forces $Q=I$, contradicting \(f\notin Q\). Thus \(\mu_R(I)=n+1\), so \(r=1\). Note that here $\tau$ is the solution of $\Phi'_{s,1,1}(n)=\HH_s'(n)-\HH_{s-1}'(n)=0$, which is thus identified with $\frac{n+1}{2}$ by \cite[Lemma 3.9]{lower-bound-s}. The result then follows by Theorem \ref{thm-nsubset-s-closure}.
\end{proof}

Consequently, following the results in \cite[Theorem 1.4, Theorem 3.5, Proposition 3.11]{mu-tau}, we have
\begin{corollary}\label{cor-mu-tau-p}
    For an isolated hypersurface singularity $f$ of dimension $n-1$, $n \geq 3$ and characteristic $p>0$, $$\frac{\mu(\mathcal{O}_f)}{\tau(f)}<\frac{1}{\mathcal{H}_{\frac{n+1}{2}}(n)-\mathcal{H}_{\frac{n-1}{2}}(n)},$$ where $\mu(\mathcal{O}_f)$ is the generalized Milnor number, $\tau(f)$ is the Tjurina number, see \cite{mu-tau}.
\end{corollary}

Together with the equality examples in dimension one, this completely answers all questions concerning the bound of the quotient of the Milnor number to the Tjurina number for hypersurface singularities in arbitrary characteristic.

Furthermore, we can get the following corollary using an argument similar to that in \cite[Proposition 3.11]{mu-tau}. It replaces the original $\leq$ by $<$ in \cite[Proposition 5.7]{Lech} (or \cite[Theorem 2.10]{lower-bound-s}, \cite[Theorem 3.1]{e_hk-eI-equality}) in the regular case. 

\begin{corollary}\label{cor-colength-e-p}
    Let $(R,\m,k)$ be a regular Noetherian local ring of dimension $n \geq 3$ and characteristic $p>0$. Let $I$ be an $\m$-primary ideal of $R$ such that $r=\mu_R(I)-n \geq 1.$ Let $1<\tau \leq \frac{n+1}{2}$ be the unique solution of $\Phi'_{s,1,r}(n)=\HH_s'(n)-r\HH_{s-1}'(n)=0$ as defined in Corollary \ref{cor-hb-ha-qua}. Then we have $$\lambda_R(R/I)=e_{HK}(I)>e(I)(\HH_\tau(n)-r\HH_{\tau-1}(n)) \geq e(I)(\HH_s(n)-r\HH_{s-1}(n)) $$ for $s \geq 1$.
\end{corollary}
\begin{proof}
    By Cohen’s structure theorem, we may write $\hat{R}=k[[x_1,\dots,x_n]]$. Let $\overline{k}$ be an algebraic closure of $k$, and set $S=\overline{k}[[x_1,\dots,x_n]]$. Then $R \rightarrow S$ is faithfully flat, and the maximal ideal of $S$ is $\m S$, see \cite[Remark 3.15]{Meng}. Consequently,
    \[
\lambda_S(S/IS)=\lambda_R(R/I),\qquad
e_S(IS)=e_R(I),\qquad
\mu_S(IS)=\mu_R(I).
\]
Therefore, we may assume that \(R\) is complete and \(F\)-finite with infinite residue field.  Thus $I$ has a reduction $Q$ that is a parameter ideal.
    
    $\lambda_R(R/I)$ coincides with $e_{HK}(I)$ when $R$ is regular. \cite[Lemma 3.9]{lower-bound-s} has shown that $\tau$ is a maximum point. Hence, we only need to show $e_{HK}(I)>e(I)(\HH_\tau(n)-r\HH_{\tau-1}(n))$, while by \cite[Theorem 2.10]{lower-bound-s} $e_{HK}(I)\geq e(I)(\HH_\tau(n)-r\HH_{\tau-1}(n))$ already holds. Assume that $e_{HK}(I)=e(I)(\HH_\tau(n)-r\HH_{\tau-1}(n))$. By \cite[Theorem 2.10]{lower-bound-s} we also know that $h_\tau(Q,I)\geq e(I)(\HH_\tau(n)-r\HH_{\tau-1}(n))$, which implies $e_{HK}(I)=h_\tau(Q,I)$.
    Therefore, since every ideal of a regular local ring is tightly closed, \cite[Theorem 7.9]{Meng} yields $$c^I(Q)=\min\{t>0:h_t(Q,I)=e_{HK}(I)\} \leq \tau.$$
    Since \(Q\) is a reduction of \(I\), we have \(c^I(Q)=c^I(I)\), contradicting Theorem \ref{thm-cJI} and Proposition \ref{prop-condition-hold}.
\end{proof}

\begin{remark}
    The condition $n \geq 3$ is necessary. If $n=2$, take $R=k[[x,y]],I=(x^{2p},y^{2p},x^py^p)$, where $\chara k=p>0$. Then $\lambda_R(R/I)=3p^2,\ e(I)=4p^2$ and $\lambda_R(R/I)=e(I)(\HH_{\frac{3}{2}}(2)-\HH_{\frac{1}{2}}(2))$. 
\end{remark}

Using the standard method of reduction modulo $p$, we can easily show the corresponding result in zero characteristic. 

\begin{corollary}\label{cor-colength-e-0}
    Let $(R,\m,k)$ be a regular Noetherian local ring of dimension $n \geq 3$ and equicharacteristic $0$. Let $I$ be an $\m$-primary ideal of $R$ such that $r=\mu_R(I)-n \geq 1.$ Let $1<\tau \leq \frac{n+1}{2}$ be the unique solution of $\Phi'_{s,1,r}(n)=\HH_s'(n)-r\HH_{s-1}'(n)=0$ as defined in Corollary \ref{cor-hb-ha-qua}. Then we have $$\lambda_R(R/I)>e(I)(\HH_\tau(n)-r\HH_{\tau-1}(n)) \geq e(I)(\HH_s(n)-r\HH_{s-1}(n)) $$ for $s \geq 1$.
\end{corollary}
\begin{proof}
    Since completion preserves colength, multiplicity, and the minimal number of generators, we may assume that $R=k[[x_1,\dots,x_n]]$ by Cohen structure theorem. Choose $N$ such that $\m^N \subseteq I$. Then the ideal $J=I \cap k[x_1,\dots,x_n]$ satisfies $\m^N \subseteq JR \subseteq \m,\ JR=I$. Let $A \subseteq k$ be the finitely generated $\mathbb{Z}$-algebra, generated by all the coefficients of the generators of $J$. Then $A$ is a Noetherian domain. Set $B=A[x_1,\dots,x_n], \n=(x_1,\dots,x_n)B$ and let $J_A$ be the ideal in $B$ generated by the generators of $J$, adjoining all monomials of degree \(N\). It follows that $\n^N \subseteq J_A \subseteq \n$, hence both $B/J_A$ and $B/\n J_A$ are finite $A$-modules. By Grothendieck’s generic freeness lemma \cite[Theorem 14.4]{Eisenbud}, after replacing \(A\) by \(A[a^{-1}]\) for some \(0\ne a\in A\), we may assume that both modules are finite free over \(A\).

    Choose a maximal ideal $\p \subset A$. Since $A$ is a finitely generated $\mathbb{Z}$-algebra, the residue field $\kappa=A/\p$ is finite. Set $$R_{\p}=\kappa[[x_1,\dots,x_n]],\ \m_\p=(x_1,\dots,x_n)R_\p,\ I_\p=J_AR_\p.$$ The freeness given above implies $$\lambda_{R_\p}(R_\p/I_\p)=\lambda_R(R/I),\quad \lambda_{R_\p}(R_\p/\m_\p I_\p)=\lambda_R(R/\m I),$$ since these lengths are the dimensions of the corresponding fibers of \(B/J_A\) and \(B/\mathfrak nJ_A\). Therefore, $$\mu_{R_\p}(I_\p)=\mu_R(I)=n+r$$ since $\mu_R(I)=\lambda_R(R/\m I)-\lambda_R(R/I)$.

    For each $t \geq 1$, \(B/J_A^t\) is finite over \(A\). By the upper semicontinuity of fiber dimensions for coherent sheaves, its fiber dimension at \(\mathfrak p\) is at least its dimension at the generic point of \(\operatorname{Spec}A\). Hence
    $$\lambda_{R_\p}(R_\p/I_\p^t)=\dim_\kappa((B/J_A^t)\otimes_A \kappa) \geq \dim_k((B/J_A^t)\otimes_A k)=\lambda_R(R/I^t). $$ Multiplying by \(n!/t^n\) and letting \(t\to\infty\), we obtain $e(I_{\mathfrak p})\ge e(I).$ Therefore, $\lambda_R(R/I)=\lambda_{R_\p}(R_\p/I_\p)>e(I_\p)(\HH_\tau(n)-r\HH_{\tau-1}(n)) \geq e(I)(\HH_s(n)-r\HH_{s-1}(n))$ for $s \geq 1$, where the first inequality follows from Corollary \ref{cor-colength-e-p}.

\end{proof}

\subsection{Test ideals, multiplier ideals, and singularity thresholds}
We now express the threshold bounds in terms of test ideals and multiplier ideals, and then derive inequalities involving the \(F\)-pure and log canonical thresholds. We first work over an \(F\)-finite regular local ring of characteristic \(p>0\).

\begin{definition}\textup{(\cite[Definition 2.9]{test-ideal-fpt})}
    Let \((R,\mathfrak m,k)\) be an $F$-finite regular local ring of dimension \(n\) and characteristic \(p>0\). Let \(I\subseteq R\) be a nonzero ideal. For \(t\geq 0\), the test ideal of the pair \((R,I^t)\) is
\[
\tau(I^t)
=
\bigcup_{e\ge1}
(I^{\lceil tp^e\rceil})^{[1/p^e]},
\]
where \(K^{[1/p^e]}\) denotes the smallest ideal \(L\subseteq R\) such that \(K\subseteq L^{[p^e]}\). The union stabilizes for sufficiently large \(e\).
\end{definition}

We have the following relation between test ideals and $F$-thresholds.
\begin{theorem}\label{thm-testideal-fthreshold}\textup{(\cite[Proposition 2.29]{test-ideal-fpt})}
    With the notation above, let \(J\subsetneq R\) satisfy \(I\subseteq\sqrt J\). Then
\[
\tau(I^t)\subseteq J
\quad\Longleftrightarrow\quad
t\ge c^J(I).
\]
Equivalently,
\[
c^J(I)
=
\min\{
t\ge0 : \tau(I^t)\subseteq J
\}.
\]
\end{theorem}

Combining Theorem \ref{thm-testideal-fthreshold} with Corollary \ref{cor-cII}, we directly get
\begin{corollary}\label{cor-test-ideal-cII}
    Assume $\dim R=n \geq 3$. Let $I$ be an $\m$-primary ideal of $R$. If $\mu_R(I) \leq n+1$, then $$\tau(I^{\frac{n+1}{2}}) \nsubseteq I.$$
\end{corollary}

\begin{definition}\textup{(\cite[Definition 9.2.3]{multiplier-ideal})}
    Let \(X\) be a smooth variety over a field \(k\) of characteristic zero, and let \(I\subseteq\mathcal O_X\) be a nonzero coherent ideal sheaf. Choose a log resolution \(\pi:Y\to X\) of \(I\), and write
$I\mathcal O_Y=\mathcal O_Y(-D).$
For \(t\ge0\), the multiplier ideal of \(I\) with exponent \(t\) is
\[
\mathcal J(I^t)
=
\pi_*\mathcal O_Y\bigl(K_{Y/X}-\lfloor tD\rfloor\bigr),
\]
where \(K_{Y/X}=K_Y-\pi^*K_X\).
\end{definition}

Test ideals are related to multiplier ideals through reduction modulo \(p\): for each fixed exponent \(t\ge0\), the reduction of the multiplier ideal \(\mathcal J(I^t)\) agrees with the test ideal \(\tau(I_p^t)\) after a suitable localization of the coefficient ring, see \cite[Theorem 3.2]{test-multiplier}. Applying this comparison to Corollary \ref{cor-test-ideal-cII} yields the following corollary.

\begin{corollary}
    Let \((R,\m)\) be a regular local ring of dimension \(n\ge3\), essentially of finite type over a field of characteristic zero. Let \(I\) be an \(\m\)-primary ideal of $R$. If $\mu_R(I) \leq n+1$, then \[
\mathcal J(I^{\frac{n+1}{2}})\not\subseteq I.
\]
\end{corollary}

\begin{definition}\textup{(\cite[Proposition 2.7]{test-multiplier}, \cite[Definition 1.16]{lct})}
    Let \((R,\mathfrak m)\) be a regular local ring, and let \(I\subset R\) be a nonzero proper ideal.

If \(R\) is \(F\)-finite of characteristic \(p>0\), the \(F\)-pure threshold of \(I\) is defined as
\[
\operatorname{fpt}(I)
=
\sup\{t\ge 0\mid \tau(I^t)=R\}
=
\min\{t\ge 0\mid \tau(I^t)\ne R\}.
\]
This agrees with the \(F\)-threshold \(c^{\mathfrak m}(I)\).

If \(R\) is essentially of finite type over a field of characteristic zero, the log canonical threshold of \(I\) is defined in terms of multiplier ideals by
\[
\operatorname{lct}(I)
=
\sup\{t\ge 0\mid \mathcal J(I^t)=R\}
=
\min\{t\ge 0\mid \mathcal J(I^t)\ne R\}.
\]
\end{definition}

Combining the strict colength inequality with the known multiplicity bound in terms of the \(F\)-pure threshold yields the following estimate.

\begin{corollary}\label{cor-colength-fpt}
    Let \((R,\mathfrak m)\) be an \(F\)-finite regular local ring of characteristic \(p>0\) and dimension \(n\ge3\). Then every \(\mathfrak m\)-primary ideal \(I\) satisfies
    $$\lambda_R(R/I)\mathrm{fpt}(I)^n > \frac{n^n}{n!(1-\mu_R(I)^{-\frac{1}{n-1}})^{n-1}}.$$
\end{corollary}
\begin{proof}
    Set $r=\mu_R(I)-n$. If $r \geq 1$, by Corollary \ref{cor-colength-e-p}, 
    \begin{equation}
        \begin{aligned}
            \lambda_R(R/I)&> e(I)\max_{s \geq 1}\{\HH_s(n)-r\HH_{s-1}(n)\}\\
            & \geq e(I)\max_{s \geq 1}\{\frac{s^n-\mu_R(I)(s-1)^n}{n!}\}\\
            &=\frac{e(I)}{n!(1-\mu_R(I)^{-\frac{1}{n-1}})^{n-1}}.
        \end{aligned}
    \end{equation}
    The second inequality is from Lemma \ref{lemma-H_s}\ref{item-H-scaling}. On the other hand, we have $$e(I)\mathrm{fpt}(I)^n \geq n^n $$ by \cite[Proposition 4.5]{e-fpt}. Then we get the desired result.

    If $r=0$, then $I$ is a parameter ideal. Then $e(I)=\lambda_R(R/I)$. The result follows easily.
\end{proof}

Specifically, we get
\begin{corollary}\label{cor-F-regular}
    With the notation above, if $n=\dim R=3,\ \mu_R(I) \leq 4$, then $$\lambda_R(R/I)\mathrm{fpt}(I)^3>18.$$ In particular, if the pair $(R,I^1)$ is not strongly $F$-regular, then $\lambda_R(R/I) \geq 19$. 
\end{corollary}
\begin{proof}
    The first assertion follows directly from Corollary \ref{cor-colength-fpt}. The second one follows from \cite[Proposition 2.2(4),(5)]{e-fpt}.
\end{proof}

\begin{remark}
    The bound $19$ is sharp. Take $I=(x^3,y^3,z^3,xyz) \subseteq \mathbb{F}_p[[x,y,z]]$. Among the $3^3$ monomials with all exponents below $3$, exactly $2^3$ are multiples of $xyz$, so $\lambda_R(R/I)=19$. For every $q=p^e$, the element $(xyz)^{q-1}$ lies in $I^{q-1}\setminus\m^{[q]}$, whereas $I^q\subseteq\m^{3q}\subseteq\m^{[q]}$. Thus $\nu_I^\m(q)=q-1$ and $\mathrm{fpt}(I)=1$. This example shows that it is necessary to replace the non-strict inequality with a strict inequality in Corollary \ref{cor-colength-e-p}.
\end{remark}

Similarly we can get an analogue in characteristic zero combining \cite[Theorem 0.1]{e-lct}.
\begin{corollary}\label{cor-klt}
    Let \((R,\mathfrak m)\) be a regular local ring essentially of finite type over a field of characteristic zero and dimension \(n\ge3\). Then every \(\mathfrak m\)-primary ideal \(I\) satisfies
    $$\lambda_R(R/I)\mathrm{lct}(I)^n > \frac{n^n}{n!(1-\mu_R(I)^{-\frac{1}{n-1}})^{n-1}}.$$ If $n=\dim R=3,\ \mu_R(I) \leq 4$, then we have $$\lambda_R(R/I)\mathrm{lct}(I)^3>18.$$ In particular, if the pair $(R,I^1)$ is not Kawamata log terminal, then $\lambda_R(R/I) \geq 19$.
\end{corollary}

\section*{Acknowledgments}
AI tools were used to assist with developing proofs, performing computations, and improving
the exposition. All results have been verified by the authors. The authors thank Dr. Cheng
Meng for stimulating discussions. Zuo is supported by NSFC Grant 12671056 and BJNSF Grant
1252009.

\newcommand{\etalchar}[1]{$^{#1}$}

\end{document}